\documentclass[12pt]{article}
\usepackage{tikz}
\usepackage{amsmath,amssymb}
\usepackage{tikz}
\usepackage{color}
\usepackage{cite}
\newtheorem{theorem}{Theorem}[section]
\newtheorem{lemma}[theorem]{Lemma}
\newtheorem{corollary}[theorem]{Corollary}
\newtheorem{proposition}[theorem]{Proposition}

\newcommand{\proof}{\noindent{\bf Proof.\ }}
\newcommand{\qed}{\hfill $\square$ \bigskip}

\newcommand{\cp}{\,\square\,}

\DeclareMathOperator {\gp} {gp}

\let\deg\relax
\DeclareMathOperator {\deg} {deg}
\DeclareMathOperator {\Pl} {Pl}

\begin{document}

\title{The general position achievement game played on graphs}

\author{
Sandi Klav\v zar\thanks{Email: \texttt{sandi.klavzar@fmf.uni-lj.si}} $\ ^{a,b,c}$
\and Neethu P. K.\thanks{Email: \texttt{p.kneethu.pk@gmail.com}} $\ ^{d}$
\and Ullas Chandran S. V.\thanks{Email: \texttt{svuc.math@gmail.com}} $\ ^{d}$
}
\maketitle

\begin{center}
$^a$ Faculty of Mathematics and Physics, University of Ljubljana, Slovenia\\
\medskip
	
$^b$ Institute of Mathematics, Physics and Mechanics, Ljubljana, Slovenia\\
\medskip
	
$^c$ Faculty of Natural Sciences and Mathematics, University of Maribor, Slovenia\\
\medskip

$^d$ Department of Mathematics, Mahatma Gandhi College, University of Kerala, Thiruvananthapuram-695004, Kerala, India
\end{center}

\begin{abstract}
A general position set of a graph $G$ is a set of vertices $S$ in $G$ such that no three vertices from $S$ lie on a common shortest path. In this paper we introduce and study the general position achievement game. The game is played on a graph $G$ by players A and B who alternatively pick vertices of $G$. A selection of a vertex is legal if has not been selected before and the set of vertices selected so far forms a general position set of $G$. The player who selects the last vertex wins the game. Playable vertices at each step of the game are described, and sufficient conditions for each of the players to win is given. The game is studied on Cartesian and lexicographic products. Among other results it is proved that A wins the game on $K_n\cp K_m$ if and only if both $n$ and $m$ are odd, and that B wins the game on $G\circ K_n$ if and only if either B wins on $G$ or $n$ is even.
\end{abstract}


\medskip\noindent

{\bf AMS Subj.\ Class:} 05C12; 05C69

%%%%%%%%%%%%%%%%%%%%%%
\section{Introduction}
\label{sec:intro}
%%%%%%%%%%%%%%%%%%%%%%

The general position problem for graphs was independently introduced and researched in~\cite{manuel-2018a, ullas-2016}, but should be noted that in the case of hypercubes, it has been studied much earlier by K\"orner~\cite{korner-1995}.  Among motives for introducing the problem is the more than a century old no-three-in-line problem of Dudeney~\cite{dudeney-1917}, see also~\cite{ku-2018, misiak-2016, skotnica-2019}. For the related general position subset selection problem in computational geometry see~\cite{froese-2017, payne-2013}. 

A {\em general position set} of a graph $G = (V(G), E(G))$ is a set of vertices $S\subseteq V(G)$ such that no three vertices from $S$ lie on a common shortest path of $G$. The general position problem asks for the largest possible size of a general position set of $G$; this number is denoted by $\gp(G)$. Immediately after its introduction, the concept received a great response~\cite{bijo-2019, ghorbani-2021, klavzar-2021a, klavzar-2021b, klavzar-2019, manuel-2018b, patkos-2020, thomas-2020, tian-2021a, tian-2021b, yao-2021+}. Furthermore, in~\cite{klavzar-2021c} general position sets have been extended to general $d$-position sets, while in~\cite{klavzar-2021d} the Steiner general position problem was studied. 

In this paper we study the achievement game associated with general position sets. Achievement games have already been studied in different contexts. For instance, in a finite group two players in turn select previously unselected elements of the group, and the player who is the first to achieve a generating set from the jointly selected elements wins the game~\cite{anderson-1987, benesh-2019}. Similarly, and closer to our game, two players in turn select vertices of a finite graph, and the player who first plays such a vertex that the union of the intervals between the vertices played contains all the vertices wins the game~\cite{buckley-1987, haynes-2003, nec-1988}. 

Let $G$ be a graph. Then the {\em general position achievement game} ({\em gp achievement game} for short) is played by two players, {\em player A} and {\em player B}. The first player A chooses a vertex $v_1$. The second player B then chooses a vertex $v_2\neq v_1$. Next A picks a vertex $v_3\in V(G)\setminus \{v_1,v_2\}$ such that the set $\{v_1,v_2,v_3\}$ is a general position set in $G$. The game then proceeds along the same way and ends when there is no more vertex to be played, that is, there exists no vertex such that the general position set consisting of the already played vertices could be enlarged. The player who has played the last vertex wins the game. 

We proceed as follows. In the rest of this section, additional definitions and notation needed are recalled. In the next section we give some general results and provide several examples. Among other results we observe that A wins the gp achievement game on a bipartite graph $G$ if and only if the number of isolated vertices in $G$ is odd. In Section~\ref{sec:Cartesian} we study the game on Cartesian products, while in Section~\ref{sec:Lexico} we prove that B wins the game on the lexicographic product $G\circ K_n$ if and only if either B wins on $G$ or $n$ is even. At the end several concluding remarks are given, among which a closely related avoidance game is commented. 

All graphs considered are finite, simple, and without loops or multiple edges.  The {\em distance} $d_G(u,v)$ between vertices $u$ and $v$ of $G$ is the length of a shortest $u,v$-path. A $u,v$-path of minimum length is also called an $u,v$-{\it geodesic}.  The {\em interval} $I_G[u,v]$ between $u$ and $v$  is the set of vertices that lie on some $u,v$-geodesic of $G$. For $S\subseteq V(G)$, we set $I_G[S]=\bigcup_{_{u,v\in S}}I_G[u,v]$. A subgraph $H$ of a graph $G$ is {\em convex}  if for every $u,v\in V(G)$, every $u,v$-geodesic in $G$ lies completely in $H$.

%%%%%%%%%%%%%%%%%%%%%%
\section{Some general results and examples}
\label{sec:preliminary}
%%%%%%%%%%%%%%%%%%%%%%

The sequence of vertices played in the achievement game on a graph $G$ will be denoted by $a_1, b_1, a_2, b_2, \ldots$, that is, the vertices played by A are $a_1, a_2, \ldots$, and the vertices played by B are $b_1, b_2, \ldots$ For instance, we may say that A starts the game by playing $a_1 = x$, where $x\in V(G)$. Suppose that $x_1, \ldots, x_j$ are vertices played so far on the graph $G$. Then we say that $y\in V(G)$ is a {\em playable vertex} if $y\notin \{x_1, \ldots, x_j\}$ and $\{x_1, \ldots, x_j\} \cup \{y\}$ is a general position set of $G$. Let $\Pl_G(x_1, \ldots, x_j)$ be the set of all playable vertices after the vertices $x_1, \ldots, x_j$ have already been played; we may sometimes simplify the notation $\Pl_G(x_1, \ldots, x_j)$ to $\Pl_G(\ldots x_j)$. For instance, if $x$ and $y$ are arbitrary vertices of a path $P$, then $\Pl_P(x) = V(P)\setminus \{x\}$ and $\Pl_P(x,y) = \emptyset$. Denoting by $S$ the set of vertices $\{x_1, \ldots, x_j\}$  played so far, we may also write $\Pl_G(S)$ for $\Pl_G(x_1, \ldots, x_j)$. In the sequel we will implicitly but frequently use the following description of playable vertices. 

\begin{lemma}
\label{lem:playable}
Let $S$ be the sequence of played vertices so far in a gp achievement game on a graph $G$. Then $x\in \Pl_G(S)$ if and only if the following two conditions hold:  
\begin{enumerate}
\item[(i)] if $u,v\in S$, then $x\notin I[u,v]$, and 
\item[(ii)] if $u\in S$, then $I[x,u]\cap S = \{u\}$. 
\end{enumerate}
\end{lemma}

\proof
By definition, $x\in \Pl_G(S)$ if and only if $S\cup \{x\}$ is a general position set. Since $S$ is a general position set by the assumption of the game, deciding whether $S\cup \{x\}$ is a general position set reduces to checking the conditions (i) and (ii). 
\qed

Let us next look at some examples. Since in a complete graph every vertex subset is a general position set, A wins the gp achievement game on the complete graph $K_n$ if and only if $n$ is odd. In the course of the gp achievement game on a graph $G$ of order at least $2$, at least two vertices will be played. Hence B wins the game on graphs $G$ with $\gp(G) = 2$. As proved in~\cite{ullas-2016}, the only graphs with  $\gp(G) = 2$ are paths and $C_4$. On the other hand, the class of graphs $G$ with  $\gp(G) = 3$ has not yet been characterized. If $\gp(G) = 3$, then gp achievement game will take either two or three moves. In fact, if $\gp(G) = 3$ then B wins the gp achievement game if only only if every vertex of $G$ lies in a maximal general position set  of order $2$. Applying this observation to cycles we infer that B wins the gp achievement game on the cycle $C_{n}$, $n\ge 3$, if and only if $n$ is even. 

The following result is simple but at the same time quite useful. 

\begin{theorem}
\label{thm:up-to-gp-set} 
Let $G$ be a graph. Then the following holds. 

(i) If A has a strategy such that after the vertex $a_k$, $k\ge 1$, is played, the set $\Pl_G(\ldots a_k)\cup \{a_1, b_1, \ldots, a_k\}$ is a general position set and $|\Pl_G(\ldots a_k)|$ is even, then A wins the gp achievement game. 

(ii) If B has a strategy such that after the vertex $b_k$, $k\ge 1$, is played, the set $\Pl_G(\ldots b_k)\cup \{a_1, b_1, \ldots, b_k\}$ is a general position set and $|\Pl_G(\ldots b_k)|$ is even, then B wins the gp achievement game. 
\end{theorem}

\proof
(i) Suppose that A has a strategy such that after A plays $a_k$, the set $\Pl_G(\ldots a_k)\cup \{a_1, b_1, \ldots, a_k\}$ is a general position set. By definition, in the rest of the game only vertices from $\Pl_G(\ldots a_k)$ are playable. Moreover, each of these vertices will actually be played because $\Pl_G(\ldots a_k)\cup \{a_1, b_1, \ldots, a_k\}$ is a general position set. Just after this will be done, the game will  be finished. Since $|\Pl_G(\ldots a_k)|$ is assumed to be even, this means that A will be the last player to select a vertex. 

(ii) Follows by a parallel argument. 
\qed

For the first application of Theorem~\ref{thm:up-to-gp-set} consider the Petersen graph $P$. Suppose that after A plays some vertex $a_1$ of $P$, B plays a vertex $b_1$ adjacent to $a_1$. Then  $\Pl_P(a_1, b_1)$ consists of four vertices which, together with $a_1$ and $b_1$, form a general position set of $P$. Hence Theorem~\ref{thm:up-to-gp-set}(ii) implies that B wins the gp achievement game on the Petersen graph. As another application of Theorem~\ref{thm:up-to-gp-set} we have the following result. 

\begin{proposition}
\label{prp:complete-multi}
Let $G$ be the complete multipartite graph $K_{n_1,\ldots, n_k}$, where $k\ge 2$ and $n_i\ge 2$ for $i\in [k]$. Then A wins the gp achievement game on $G$ if and only if $k$ is odd and at least one $n_i$ is odd.  
\end{proposition}

\proof
Suppose first that $n_i$ is even for all $i\in [k]$. Let $X$ be the partition set of $G$ in which the first move $a_1$ has been played by A. Then B replies by playing a vertex $b_1\ne a_1$ from $X$. Note that $\Pl_G(a_1, b_1) = X\setminus \{a_1, b_1\}$. Since $X$ is a general position set of $G$, Theorem~\ref{thm:up-to-gp-set}(ii) applies and B wins the gp achievement game on $G$. 

Hence, the only possibility for A to win the game is that at least one $n_i$ is odd and that the first move $a_1$ is from an odd partition set $X$. Now, if B would reply by playing a vertex in $X$, then by the argument of the previous paragraph and with Theorem~\ref{thm:up-to-gp-set}(i) in hand, A would win. So it is better for B to play a vertex $b_1$ which lies in a partition set $Y\ne X$. Since $\Pl_G(a_1,b_1) = V(G)\setminus (X\cup Y)$,  the vertex $a_2$ must lie in a partition set $Z$ different from both $X$ and $Y$. Continuing in this manner, each of the subsequent played vertices belongs to its private partition set. In conclusion, if some $n_i$ is odd, then A will win if and only if $k$ is odd.
\qed

In view of Theorem~\ref{thm:up-to-gp-set}(ii) we easily infer that B wins the gp achievement game on an arbitrary connected, bipartite graph of order at least two. This observation generalizes to arbitrary bipartite graphs as stated in the next theorem, for which we need the following fact that was observed for the first time in the proof of~\cite[Theorem 5.1]{bijo-2019}.  

\begin{lemma}
\label{lem:bipartite} 
Let $G$ be a connected, bipartite graph. If $S$ is a general position set of $G$ with $|S|\geq 3$, then $S$ is an independent set.
\end{lemma} 

\begin{theorem}
\label{thm:bipartite} 
Let $G$ be a bipartite graph. Then A wins the gp achievement game on $G$ if and only if the number of isolated vertices in $G$ is odd.
\end{theorem}

\proof 
Let $k$ be the number of isolated vertices of $G$. 

First suppose that $k\ge 0$ is even and consider the following strategy of B. Whenever A selects a vertex $v$ in some component of $G$ of order at least $2$, B replies with a move on a neighbor of $v$. And whenever A plays an isolated vertex, B replies by playing another isolated vertex. Note that after two adjacent vertices of a component $H$ of $G$ are played, Lemma~\ref{lem:bipartite} implies that no additional vertex from $H$ will be played in the rest of the game. Moreover, since $k$ is even, whenever A plays an isolated vertex, there exists at least one isolated vertex which was not played yet, hence B can follow the described strategy. It follows that the game will finish when all the isolated vertices and precisely two (adjacent) vertices from each component will be played. So the number of played vertices will be even, hence B wins the gp achievement game. 

Second, let $k\ge 1$ be odd. Then A has the following strategy to win the gp achievement game. The first vertex played will be an isolated vertex. After that, the strategy of A is just as the described strategy of B in the above paragraph: whenever B plays a vertex $v$ in some component of $G$ of order at least $2$, A replies with a move on a neighbor of $v$, and if B plays an isolated vertex, A replies by playing another isolated vertex. Using parallel arguments as above, the total number of vertices played will be odd, which in turn implies that A wins the game. 
\qed

%%%%%%%%%%%%%%%%%%%%%%
\section{The game played on Cartesian products}
\label{sec:Cartesian}
%%%%%%%%%%%%%%%%%%%%%%

 The {\em Cartesian product} $G\cp H$ of graphs $G$ and $H$ has the vertex set $V(G)\times V(H)$, the vertices $(g_1,h_1),(g_2,h_2)$ being adjacent in $G\cp H$ if either $g_1g_2\in E(G)$ and $h_1=h_2$, or $g_1=g_2$ and $h_1h_2\in E(H)$. If $g\in V(G)$, then the subgraph of $G\cp H$ induced by the vertex set $\{(g,h)$ $|$ $h\in v(H)\}$ is an {\em $H$-layer} ${ }^{g}H$.  {\em $G$-layers} $G^h$ are defined analogously. If $S\subseteq V(G\cp H)$, then the {\em projection} $\pi_G(S)$ of $S$ on $G$ is the set $\{g\in V(G):\ (g,h)\in \ {\rm for\ some}\ h\in V(H)\}$. The projection $\pi_H(S)$ of $S$ on $H$ is defined analogously. 

Throughout this section we will use the following basic fact about the distance function in the Cartesian product. If $G$ and $H$ are connected graphs and $(g,h), (g',h')\in V(G\cp H)$, then the {\em distance formula} holds: 
\begin{equation}
\label{eq:distance-formula}
d_{G\cp H}((g,h), (g', h')) = d_G(g,g') + d_H(h,h')\,.
\end{equation}
Moreover, if $P$ is a $(g,h),(g',h')$-geodesic in $G\cp H$, then $\pi_G(P)$ induces a $g,g'$-geodesic in $G$ and $\pi_H(P)$ induces a $h,h'$-geodesic in $H$. The distance formula~\eqref{eq:distance-formula} implies that 
\begin{equation}
\label{eq:intervals-in-CP}
I_{G\cp H}[(g,h), (g', h')] = I_G[g,g'] \times I_H[h,h']\,.
\end{equation} 
For these results and more on the Cartesian product operation see the standard book on product graphs~\cite{hik-2011}. We will also need the following known result. 

\begin{lemma}
	{\rm\cite[Lemma 2.4]{tian-2021b}}
	\label{lemma:2.4}
Let $G$ and $H$ be connected graphs and let $R$ be a general position set of $G\cp H$. If $u=(g,h)\in R$, then $V(^{g}{H})\cap R=\{u\}$ or $V(G^h)\cap R=\{u\}$. 
\end{lemma}

We next prove two additional lemmas on general position sets in Cartesian products. 

\begin{lemma}
\label{lemma:3.7}
Let $G$ and $H$ be connected graphs and let $R\subseteq V(G\cp H)$ has the following two properties. 
\begin{enumerate}
\item[(i)]  If $(g,h)\in R$, then $V(^{g}{H})\cap R=\{(g,h)\}$ or $V(G^h)\cap R= \{(g,h)\}$.
\item[(ii)] $\pi _G(R)$ and $\pi_H(R)$ are general position sets of $G$ and $H$, respectively.
\end{enumerate}
Then $R$ is a general position set of $G\cp H$.
\end{lemma}

\proof 
Suppose on the contrary that $R$ contains three vertices $x_1 = (u_1,v_1)$, $x_2 = (u_2,v_2)$, and $x_3 = (u_3,v_3)$ such that $x_2\in I_{G\cp H}[x_1, x_3]$. Applying the distance formula~\eqref{eq:distance-formula} and the triangle inequality we can estimate as follows: 
\begin{align*}
d_{G\cp H}(x_1,x_3) &  = d_{G\cp H}(x_1, x_2) + d_{G\cp H}( x_2, x_3) \\ 
& = (d_G(u_1,u_2)+d_H(v_1,v_2)) + (d_G(u_2,u_3)+d_H(v_2,v_3)) \\
& = (d_G(u_1,u_2) + d_G(u_2,u_3)) + (d_H(v_1,v_2) + d_H(v_2,v_3)) \\
& \ge d_G(u_1,u_3)+ d_H(v_1,v_3) \\
& = d_{G\cp H}(x_1,x_3)\,.
\end{align*} 
It follows that $d_G(u_1,u_3) = d_G(u_1,u_2) + d_G(u_2,u_3)$ and $d_H(v_1,v_3) = d_H(v_1,v_2) + d_H(v_2,v_3)$. 

Suppose first that $x_1$ and $x_3$ lie in a common $G$-layer or in a common $H$-layer. By the commutativity of the Cartesian product we may without loss of generality assume that they lie in a common $H$-layer, that is, $u_1=u_3$.  Since $^{u_1}H$ is a convex subgraph of $G\cp H$ (see~\cite{hik-2011} again), it follows that $u_1 = u_2 = u_3$. Hence the vertices $v_1, v_2, v_3$ are pairwise different, and so the fact $d_H(v_1,v_3) = d_H(v_1,v_2) + d_H(v_2,v_3)$ yields a contradiction with the assumption that $\pi_H(R)$ is a general position set of $H$.

Assume second that $x_1$ and $x_3$ lie neither in a common $G$-layer nor in a common $H$-layer. Then $u_1\ne u_3$ and $v_1\ne v_3$. Assumption (i) then implies that $(u_2,v_2) \notin \{(u_1,v_3), (u_3,v_1)\}$. As a consequence, at least one of the sets $\{u_1,u_2,u_3\}$ and $\{v_1,v_2,v_3\}$ is of cardinality $3$. But then we have a contradiction for one of these sets just as in the previous paragraph. 
\qed

The converse of Lemma~\ref{lemma:3.7} does not hold. As an example consider the path $P_3$ on vertices $1, 2, 3$, and the Cartesian product $P_3\cp P_3$. Then $\{(1,2), (2,1), (2,3), (3,2)\}$ is a general position set of $P_3\cp P_3$, but neither its projection onto the first factor not the projection onto the second factor is a general position set.  

If each $H$-layer contains at most one vertex from $R$, then the conditions of Lemma~\ref{lemma:3.7} simplify as follows. 

\begin{lemma}
\label{lemma:3.8}
Let $G$ and $H$ be connected graphs and let $R\subseteq V(G\cp H)$. If $\pi _G(R)$ is a general position set in $G$ and $\pi _G(R)=|R|$, then $R$ is a general position set of $G\cp H$.
\end{lemma}

The proof of Lemma~\ref{lemma:3.8} proceeds along the same lines as the proof of Lemma~\ref{lemma:3.7} and is hence omitted. That the converse of Lemma~\ref{lemma:3.8} again does not hold, consider again  the Cartesian product $P_3\cp P_3$. Then $\{(1,1), (2,2), (3,1)\}$ is a general position set of $P_3\cp P_3$ with exactly one vertex in each of the layers with respect to the first factor, but its projection onto the first factor is not a general position set.  

\begin{lemma}
\label{lemma:3.10}
Let $G$ and $H$ be connected graphs. If for every $u\in V(G)$ there exists a vertex $v$ such that $\Pl_{G}(u,v)\cup\{u,v\}$ is a clique of even order, then B wins the gp achievement game on $G\cp H$. 
\end{lemma}

\proof 
Consider the gp achievement game on $G\cp H$. Let $a_1=(u_1,v_1)$. Then there exists a vertex $u_2\in V(G)$ such that  $\Pl_{G}(u_1,u_2)\cup \{u_1,u_2\}$ is a clique of even order. The initial strategy of B is to play $b_1=(u_2,v_1)$. Suppose that A next plays $a_2 = (u_3,v_3)$. By Lemma~\ref{lemma:2.4}, $\Pl_{G\cp H}(a_1, b_1)\subseteq V(G\cp H)\backslash (^{u_1}H\cup{ }^{u_2}H)$, hence $u_3\notin \{u_1, u_2\}$. We claim that $\{u_1,u_2,u_3\}$ is a general position set of $G$. If not, then, since $u_1u_2\in V(G)$, we may without loss of generality assume that $u_2\in I_{G}[u_1,u_3]$. But then the distance formula implies that $(u_2,v_1)\in I_{G\cp H}[(u_1,v_1),(u_3,v_3)]$. Hence the claim, which in turn implies that $u_3\in \Pl_{G}(u_1,u_2)$. Since $|\Pl_G(u_1,u_2)|$ is even and $\Pl_G(u_1,u_2,u_3) = \Pl_G(u_1,u_2)\backslash \{u_3\}$,  player B can continue the game by choosing the vertex $(u_4,v_3)$, where  $u_4 \in \Pl_{G}(u_1,u_2, u_3)$. By Lemma~\ref{lemma:3.8}, the set $S_4=\{(u_1,v_1), (u_2,v_1), (u_3,v_3), (u_4,v_3)\}$ is a general position set of $G\cp H$. Player B then continues this strategy and by repeatedly applying Lemma~\ref{lemma:3.8}, we can see that each set $S_n$ is a general position set of $G\cp H$. Also since B wins on $G$, at each stage of the game $|\Pl_{G\cp H}(\ldots b_k)|$ is even. Hence by Theorem~\ref{thm:up-to-gp-set}(ii), B wins on $G\cp H$. 
\qed

Since in a connected, bipartite graph, every pair of adjacent vertices is a maximal general position set, the following theorem follows directly from Lemma~\ref{lemma:3.10}.

\begin{theorem}
\label{thm:p-bipartite} 
Let $G$ be a connected graph and let $H$ be a connected bipartite graph with at least one edge. Then B wins the gp achievement game on $G\cp H$ after his first move. 
\end{theorem}

Theorem~\ref{thm:p-bipartite} should be compared with the main result from~\cite{tian-2021b} which asserts that if $T$ and $T'$ are trees, then $\gp(T\cp T') = \ell(T) + \ell(T')$, where $\ell(G)$ is the number of leaves of a graph $G$.    

We next resolve the gp achievement on Hamming graphs.

\begin{theorem}
\label{thm:3.11}
If $n, m\ge 2$, then A wins the gp achievement game on $K_n\cp K_m$ if and only if both $n$ and $m$ are odd.
\end{theorem}

\proof 
Let $V(K_n) = \{u_1,\ldots,u_n\}$,  $V(K_m) = \{v_1,\ldots,v_m\}$, and set $G = K_n\cp K_m$ for the rest of the proof. If one of $n$ and $m$ is even, then B wins the gp achievement game on $G$ by Lemma~\ref{lemma:3.10}. 

In the rest assume that both $n$ and $m$ are odd. We need to prove that in this case A wins the gp achievement game. The strategy of A is to achieve the following goal. After each move $a_i$, $i\ge 1$, we have that     
\begin{equation}
\label{eq:Hamming-graphs}
|\Pl_G(\ldots a_i) \cap V(K_n^{v_k})|\ 
{\rm is\ even\ and}\ 
|\Pl_G(\ldots a_i) \cap V(^{u_j}K_m)|\
{\rm is\ even}
\end{equation}
for all layers $K_n^{v_k}$ and all layers $^{u_j}K_m$ in which at least one vertex has already been played. 

By the vertex-transitivity of $G$ we may assume that $a_1=(u_1,v_1)$. Note that~\eqref{eq:Hamming-graphs} holds true after this move. For the first move $b_1=(u_i,v_j)$ of B we may, again using the symmetry of $G$, without loss of generality assume that $i=2$ and $j\in[2]$. 

Suppose first that $b_1 = (u_2,v_1)$. Then A selects $a_2=(u_3,v_1)$. Since by Lemma~\ref{lemma:2.4}, $\Pl_{G}(a_1, b_1,a_2)\subseteq V(G)\backslash (^{u_1}H\cup{ }^{u_2}H \cup {}^{u_3} H)$, the condition~\eqref{eq:Hamming-graphs}  is fulfilled after the move $a_2$. The next move of B must be in a new $K_m$-layer, say $b_2=(u_4,v_j)$. Then A replies by the vertex $a_3=(u_5, v_j)$. This is a legal move since $n$ is odd and because Lemma~\ref{lemma:3.8} guarantees that the so far selected vertices form a general position set of $G$.  The game then continues in this manner, that is, whenever it is B's turn, he must select a vertex $x$ from  some new $K_m$-layer, and then A replies with a playable neighbor of $x$ in the corresponding $K_n$-layer. As $n$ is odd, A will play the last vertex. 
 
Suppose second that $b_1 = (u_2,v_2)$. In this case A replies by picking $a_2=(u_3,v_3)$. Then $\Pl_G(a_1,b_1,a_2) = V(G) \setminus (\{u_1,u_2,u_3\} \times \{v_1,v_2,v_3\})$ and~\eqref{eq:Hamming-graphs} is fulfilled after the move $a_2$. In the sequel of the game, if B plays a vertex such that it is the first vertex played in the two layers in which it lies, then A replies with another such vertex. Note that this is possible as both $n$ and $m$ are odd. After each such move of A, the conditions~\eqref{eq:Hamming-graphs} remain fulfilled. Suppose now that at some point of the game, B selects a vertex in a $K_n$-layer in which at least one vertex has been played earlier. Because before this move~\eqref{eq:Hamming-graphs} holds, $A$ can reply by playing a vertex from the same  $K_n$-layer. Now, in this $K_n$-layer exactly two less vertices are playable, so the number of playable vertices in the layer is even (possibly zero). Moreover, in the two $K_m$-layers, in which the last two moves were played, no vertex is now playable, hence~\eqref{eq:Hamming-graphs} holds also for these two layers. In the case when at some point of the game, B selects a vertex in a $K_m$-layer in which at least one vertex has been played earlier, A proceeds analogously, that is, he plays next a vertex from the same  $K_m$-layer. Following this strategy, A wins the game. 
\qed

Lemma~\ref{lemma:3.10} immediately implies the following. 

\begin{corollary}
If $n$ is even and $G$ is a connected graph, then B wins the gp achievement game on $K_n\cp G$. 
\end{corollary}

On the other hand, if $n$ is odd and $G$ is a connected graph, the the outcome of the gp achievement game on $K_n\cp G$ appears more involved. This statement is in part justified by the following result. 

\begin{theorem}
If $m\ge 3$, then A wins the gp achievement game on $K_3 \cp C_m$ if and only if $m\in \{3,5\}$.
\end{theorem}

\proof 
Let $V(K_3)=\{u_1,u_2,u_3\}$ and let $V(C_m) = \{v_1, \ldots, v_{2k+1}\}$, where the edges are in natural order. Set $G = K_3 \cp C_m$ for the rest of the proof. 

If  $m$ is even, then Theorem \ref{thm:p-bipartite} implies that B wins the gp achievement game on $G$. If $m=3$, then by Theorem~\ref{thm:3.11} we know that A wins the game. Consider next the case $m=5$. We are going to prove that A wins the game by considering all possibilities (up to symmetry). Let A start the game with $a_1 = (u_1, v_1)$. Then, up to the symmetry of $G$ and having in mind that each vertex of $G$ is at distance at most $3$ from $a_1$, we need to consider the following replies of B: $(u_1,v_2)$,  $(u_1,v_3)$, $(u_2,v_1)$, $(u_2,v_2)$, and $(u_2,v_3)$. 
If $b_1 = (u_1,v_2)$, then A selects $a_2 = (u_1,v_4)$ and wins the game.    
If $b_1 = (u_1,v_3)$, then A selects $a_2 = (u_1,v_5)$ and wins the game.
If $b_1 = (u_2,v_1)$, then the move $a_2 = (u_3, v_1)$ finishes the game.  
If $b_1 = (u_2,v_2)$, then A replies by $a_2 = (u_3, v_3)$.  Then we have three subcases: if $b_2 = (u_1,v_4)$, then A plays $a_3 = (u_2,v_5)$; if $b_2 = (u_2,v_4)$, then A plays $a_3 = (u_3,v_5)$; and if $b_2 = (u_3,v_4)$, then A plays $a_3 = (u_2,v_5)$. In each of the subcases, A wins. Finally, if $b_1 = (u_2,v_3)$, then A selects $a_2 = (u_3, v_2)$. Similarly as in the second case we now see that B cannot win with the move $b_2$, while afterwards A wins with his third move. 

It remains to prove that B wins on $G$ when $m\ge 7$ is odd. By the vertex-transitivity of $G$ we may assume that A starts with the vertex $(u_1,v_1)$. Then B picks the vertex $b_1 = (u_2,v_2)$. From here on, we distinguish two cases. 

\medskip\noindent 
\textbf{Case 1}: $a_2 = (u_3,v_i)$, where $i>2$.\\
We may without loss of generality assume that $i\leq k+1$. By Lemma~\ref{lem:playable} and by~\eqref{eq:intervals-in-CP}, the move $b_2 = (u_3,v_s)$, where $s\geq k+2$, is a legal move of B. Let $a_3=(u_r,v_l)$. If $r=1$ and $l>k+2$, then $a_1\in I_G[a_3,b_1]$, and if $r=2$ and $l\leq k+1$, then $b_1\in I_G[a_1,a_3]$. If $r=3$, then either $\{a_1, a_2, a_3\}$ or $\{a_1, b_2, a_3\}$ is not a general position set of $G$. These cases imply that 
\begin{equation}
\label{eq:claim2}
\Pl_G(a_1,b_1,a_2,b_2)\subseteq (u_1 \times \{v_1,\ldots,v_{k+2}\}) \cup (u_2 \times \{v_{k+2}, \ldots, v_{2k+1}\})\,. 
\end{equation}
We can without loss of generality assume that A continues by playing $a_3 = (u_2,v_j)$, where $j>k+1$. Since $m \ge 7$, B can choose $b_3 = (u_1,v_t)$ with $2<t\leq k+2$ and $i\neq t\neq j$. Again by Lemma~\ref{lem:playable} and~\eqref{eq:intervals-in-CP}, $S=\{a_1,b_1,a_2,b_2,a_3,b_3\}$ is a general position set of $G$. We claim that $S$ is a maximal general position set. Suppose on the contrary that $a_4 = (u_g,v_h)$ is a legal move. Then by~\eqref{eq:claim2}, $g\neq 3$. If $g=1$, then applying~\eqref{eq:claim2} again, $h\leq k+2$. But then by~\eqref{eq:intervals-in-CP}, $\{b_1,b_3,a_4\}$ is not be a general position set of $G$. And if $g=2$, then again by~\eqref{eq:claim2}, $h>k+1$. If $h=k+2$ or $j=k+2$, then  $\{a_1,a_3,a_4\}$ is not be a general position set. And if $h>k+2$ and $j>k+2$, then $\{b_1,a_3,a_4\}$ is not be a general position set in $G$. This proves the claim which in turn finishes th argument for Case 1. 

\medskip\noindent
\textbf{Case 2}:  $a_2=(u_i,v_j)$, $i\in [2]$ or $j\in [2]$.\\
First suppose that $i\in[2]$, say $i=1$. Then clearly $j\neq 2$. If $j>k+2$, then $a_1\in I_G[a_2,b_1]$. Hence $2<j\leq k+2$. Set $b_2 = (u_3,v_2)$. By~Lemma \ref{lem:playable} and~\eqref{eq:intervals-in-CP}, the set $\{a_1,b_1,a_2,b_2\}$ is a general position set. Using Lemma~\ref{lemma:2.4}, we get $\Pl_G(a_1, b_1,a_2,b_2)\subseteq\, ^{u_1}C_m$. If A can choose $a_3 = (u_1,v_r)$, then, as mentioned above, $r\leq k+2$. But then by~\eqref{eq:intervals-in-CP}, $\{b_1,a_2,a_3\}$ is not be a general position set. Hence B wins the game. Suppose second that $j\in[2]$, say $j=2$. Then clearly $i=3$. Hence by Lemma~\ref{lemma:2.4}, $\Pl(a_1, b_1,a_2)\subseteq$ $ ^{u_1}C_m$. Now B can choose $b_2 = (u_1,v_s)$  with $2<s\leq k+2$. Then as in the case $i=1$, B wins the gp achievement game on $G$.
\qed

%%%%%%%%%%%%%%%%%%%%%%
\section{The game played on lexicographic products}
\label{sec:Lexico}
%%%%%%%%%%%%%%%%%%%%%%

The {\em lexicographic product} $G\circ H$ of graphs $G$ and $H$ has the vertex set $V(G)\times V(H)$, vertices $(g,h)$ and $(g',h')$ being adjacent if either $gg'\in E(G)$, or $g=g'$ and $hh'\in E(H)$. Layers and projections are defined for the lexicographic product in the same way as they are defined for the Cartesian product. The distances in lexicographic products can be computed as follows, see~\cite[Proposition 5.12]{hik-2011}.  

\begin{proposition}
\label{prop:Products:LexDist}
If $(g,h)$ and $(g',h')$ are two vertices of $G\circ H$, then
$$d_{G\circ H}\left((g,h),(g',h')\right)=
\left\{\begin{array}{ll}
d_G(g,g'); & \mbox{$g\ne g'$}\,,\\
d_H(h,h'); & g=g', \deg_G(g)=0\,,\\
\min\{d_H(h,h'), 2\}; & g=g', \deg_G(g)\ne 0\,.\\
\end{array}\right.$$
\end{proposition}

\begin{lemma}
\label{lemma:3.15}
If $G$ and $H$ are connected graphs, and $S$ is a general position set of $G\circ H$, then $\pi_G(S)$ is a general position set of $G$.
\end{lemma}

\proof 
Let $S$ be a general position set of $G\circ H$, and suppose on the contrary that $\pi_G(S)$ is not a general position set of $G$. Then there exist vertices $(u_1,v_1)$, $(u_2,v_2)$, and $(u_3,v_3)$ from $S$ such that $u_2\in I_{G}[u_1,u_3]$. Since $u_1$,$ u_2$, and $u_3$ are pairwise distinct, Proposition~\ref{prop:Products:LexDist} yields
\begin{align*}
d_{G\circ H}((u_1,v_1),(u_3,v_3)) & = d_{G}(u_1,u_3) = d_{G}(u_1,u_2)+d_{G}(u_2,u_3) \\ &= d_{G\circ H}((u_1,v_1),(u_2,v_2))+d_{G\circ H}((u_2,v_2),(u_3,v_3)),
\end{align*}
which is not possible as $S$ is a general position set of $G\circ H$. 
\qed

\begin{theorem}
\label{thm:3.17}
If $G$ is a connected graph, then B wins the gp achievement game on $G\circ K_n$ if and only if either B wins on $G$ or $n$ is even.
\end{theorem}

\proof 
We first claim that $S\subseteq V(G\circ H)$ is a maximal general position set of $G\circ K_n$ if and only if $S= S_G \times V(K_n)$, where $S_G$ is a maximal general position set of $G$. Proposition~\ref{prop:Products:LexDist} and Lemma~\ref{lemma:3.15} imply that if $S_G$ is a maximal general position set of $G$, then $S_G\times K_n$ is a maximal general position set of $G\circ K_n$. On the other hand, let $S$ be a maximal general position set of $G\circ K_n$. The maximality implies that if $(u,v)\in S$, then $V(^uK_n) \subseteq S$. Hence $S=\pi_G(S)\times V(K_n)$, and clearly, $\pi_G(S)$ is a maximal general position set in $G$. This proves the claim. 

Suppose that $n$ is odd and that A wins the game on $G$. We will show that then A wins also on $G\circ H$. The strategy of A is the following. First he selects a vertex $(u,v)$, where $u$ is an optimal start vertex in the game played on $G$ and $v$ is an arbitrary vertex of $K_n$. After that, A replies to moves of B in the following way. Whenever B selects a vertex $b_i$ from a $K_n$-layer, from which no vertex was played earlier, A replies with a vertex $a_{i+1}$ such that $\pi_G(a_{i+1})$ is an optimal reply of A to the move $\pi_G(b_{i})$ of B played in $G$. On the other hand, whenever $b_i$ belongs to a previously visited $K_n$-layer, A replies by choosing a vertex $a_{i+1}$ such that $\pi_G(a_{i+1}) = \pi_G(b_{i})$. Note that this is possible since $n$ is odd. Because A wins on $G$, the described strategy implies that vertices from an odd number of $K_n$-layers will be played during the game. In addition, the above claim implies that all the vertices from these $K_n$-layers will be played, hence in total odd number vertices will be played. We conclude that A wins the game on $G\circ H$ when $n$ is odd and A wins the game on $G$.

It remains to prove that in the other cases B has a winning strategy. If $n$ is even, then the claim implies that an even number of vertices will be played during the game. This means that B wins. And if B has a winning strategy on $G$, then B follows a similar strategy as A in the previous paragraph. Whenever A plays in some new $K_n$-layer, B replies optimally (with respect to the projected game on $G$) in a new $K_n$-layer, and whenever A plays in an already visited $K_n$-layer, B plays in some other already visited $K_n$-layer in which not all vertices has been played yet. In this way an even number of $K_n$-layer will be visited during the game and then the claim implies that B wins the game. 
\qed

%%%%%%%%%%%%%%%%%%%%%%
\section{Concluding remarks}
\label{sec:conclude}
%%%%%%%%%%%%%%%%%%%%%%

The main message of this paper is the following. If at an early stage of the gp achievement game played on $G$, one of the players has a possibility to play a vertex that significantly reduces the number of playable vertices during the rest of the game, then it is often easier to analyze the gp achievement game on $G$ than the general position number of $G$. On the other hand, we have also seen that many challenging problems concerning the gp achievement game remain open. 

In parallel to achievement games one can also consider avoidance games. In our particular case of the gp achievement game, the gp avoidance game is defined analogously, the only difference is that in the general position avoidance game the player who has played the last vertex loses the game. So the gp achievement game and the gp avoidance game are much similar, however, they are in general independent. For instance, one can check that player B wins both games on the cycle $C_6$. Hence it would be interesting to investigate also the gp avoidance game and to compare it with the gp achievement game. 

%%%%%%%%%%%%%%%%%%%%%%
\section*{Acknowledgments}
%%%%%%%%%%%%%%%%%%%%%%

Sandi Klav\v{z}ar acknowledges the financial support from the Slovenian Research Agency (research core funding P1-0297, and projects N1-0095, J1-1693, J1-2452). Neethu P K acknowledges the Council of Scientific and Industrial Research(CSIR), Govt.\ of India for providing financial assistance in the form of Junior Research Fellowship.

\baselineskip13pt

%%%%%%%%%%%%%%%%%%%%%%

\end{document}